\documentclass[12pt]{article}
\usepackage[cp1251]{inputenc}
\usepackage[english]{babel}
\usepackage[intlimits]{amsmath}
\usepackage{amsfonts}
\usepackage{amssymb}
\usepackage{graphicx}
\usepackage{amsthm}
\usepackage {amsmath}

\begin{document}

\newtheorem{theorem}{Theorem}
\newtheorem{lemma}{Lemma}
\newtheorem{Prop}{Proposition}
\newtheorem{zam}{Remark}
\newtheorem{sled}{Corollary}
\newtheorem{opr}{Definition}
\centerline{\Large \bf Continued fractions with minimal
remainders} \vspace*{4mm} \centerline{\Large \bf Jabitskaya E.\,N.
\footnote{ Research is supported by RFBR grant No 09-01-00371-a} }
\vspace*{4mm} Consider the representation of a rational number in
the form, associated with "centered" Euclidean algorithm. We prove
a new formula for the limit distribution function for sequences of
rationals with bounded sum of partial quotients.

\section{Introduction and main results}
 The classical
Euclidean algorithm for $a \in \mathbb{Z}$ and $b \in \mathbb{N}$
uses the division of the form
$$
a = bq + r, \quad q \in \mathbb{Z}, \quad b>0, \quad 0 \leqslant r <b.
$$
It leads to a continued fraction expansion of a real number
\begin{equation}\label{ordcf}
x = [ a_0; a_1, a_2, \ldots, a_m, \ldots] =
a_0+\cfrac{1}{a_1+\cfrac{1}{a_{2}
+ \ldots + \cfrac{1}{a_m + \ldots}}},
\end{equation}
where $a_0 \in \mathbb{Z}$, $a_j \in \mathbb{N}$ for $j \geqslant 1.$
Numbers $a_i$ are called partial quotients of fraction (\ref{ordcf}).

For $x \in \mathbb{Q}$ the representation (\ref{ordcf}) is finite.
We assume for the uniqueness  that the last partial quotient $a_l
 $ is greater or equal than $2$. Let
$$S^{[0]}(a/b):= a_0 + a_1 + \ldots + a_l.$$
Define the set

\begin{equation}\label{F_n}
\mathcal{F}_n := \left\{
x \in \mathbb{Q}, x \in [0,1] : S^{[0]}(x) \leqslant n +
1\right\}.
\end{equation} 
The limit distribution function
$$F^{[0]} (x) := \lim_{n \to \infty} \frac{\sharp \left\{ \xi \in \mathcal{F}_n:
\xi \leqslant x\right\}}{\sharp \mathcal{F}_n},\quad x \in [0,1]
$$
is the famous Minkowski's question mark function $?(x)$.
Properties of $?(x)$ were investigated  for example in
\cite{Denjoy}, \cite{Para}  and \cite{4}.

There are different kinds of Euclidean algorithms. For example,
"by-excess" Euclidean algorithm uses the division "by-excess"
$$
a = bq + r, \quad -b < r \leqslant 0,
$$
This algorithm leads to regular reduced continued fraction
(\cite{Fink}, \cite{Jabitskaya1}) expansion of a real number $x$,
that is
\begin{equation}\label{rrcf}
x = [[ a_0; a_1, a_2, \ldots, a_m, \ldots]] =
a_0-\cfrac{1}{a_1-\cfrac{1}{a_{2}
- \ldots - \cfrac{1}{a_m - \ldots}}},
\end{equation}
where $a_0 \in \mathbb{Z}$, $a_j \geqslant 2$ for $j \geqslant 1.$
Numbers $a_i$ are called partial quotients of fraction (\ref{rrcf}).

For a rational $x$ representation (\ref{rrcf}) is finite. For
rational $x$ we  denote the sum of partial quotients in the
representation of $x$ in the form (\ref{rrcf}) by $S^{[1]}(x)$. We
put
\begin{equation*}
\Xi_n: = \left\{x \in \mathbb{Q}, x \in [0,1] : S^{[1]}(x) \leqslant n + 2\right\}
\end{equation*}
Consider the limit distribution function
$$F^{[1]}(x): = \lim_{n \to \infty} \frac{\sharp \left\{ \xi \in \Xi_n:
\xi \leqslant x\right\}}{\sharp \Xi_n},\quad x \in [0,1],
$$

In 1995 R.\,F.\,Tichy and J.\,Uitz \cite{T-U} considered a one
parameter family $g_{\lambda}(x) $, $\lambda \in (0,1)$, $x \in
[0,1]$, of singular functions. Functions $F^{[0]}$, $F^{[1]}$
belong to this family with $\lambda = \frac 1 2$ and $\lambda =
\frac{3-\sqrt{5}}{2}$ correspondingly. Similar functions
$\kappa(x, \alpha)$, $x \in [0, \infty)$, $\lambda \in (0,1)$,
were introduced by A.\,Denjoy    \cite{Denjoy}  much more earlier
in 1938. For $x \in [0,1]$ functions $\kappa(x, \alpha)$ and
$g_{\lambda}(x)$ are related in the following way:
$$\kappa(x, \alpha) = 1 - (1 - \alpha)g_{1 - \alpha}(x).$$
In the same paper A.\,Denjoy proved that
$$\kappa(x, \alpha) = \alpha^{a_0} - \alpha^{a_0}(1 - \alpha)^{a_1} +
\alpha^{a_0 + a_2}(1 - \alpha)^{a_1 +a_3} - \ldots,$$
where $a_0, a_1, \ldots, a_m, \ldots$ are partial quotients of representation
(\ref{ordcf}) of number $x$.
Similar formula for $g_{\lambda}(x)$, $x \in (0,1)$ is given in the paper
\cite{Jabitskaya1}:
\begin{multline}\label{general_formula_salem}
g_{\lambda} (x) = {\lambda}^{a_1 - 1} -{\lambda}^{a_1 - 1} (1-\lambda)^{a_2 } +
{\lambda}^{a_1-1} (1-\lambda)^{a_2 } \lambda^ {a_3}  - \ldots
+\\+(-1)^{m+1} {\lambda}^{\sum\limits_{{1 \leqslant i \leqslant m}
\atop {i \equiv 1 \mod 2}} a_i- 1}
(1-\lambda)^{\sum\limits_{{1 \leqslant i \leqslant m}\atop
{i \equiv 0 \mod 2}} a_i} + \ldots.
\end{multline}
For $\lambda = \frac{1}{2}$ formula (\ref{general_formula_salem})
gives a well-known result by R.\,Salem \cite{Salem}, namely
$$
g_{1/2} (x) = ?(x) = \frac{1}{2^{a_1 - 1}} - \frac{1}{2^{a_1 + a_2 -1}} +
\frac{1}{2^{a_1 + a_2 + a_3 - 1}}- \ldots.
$$

Let us consider a "centered" version of the Euclidean algorithm.
This algorithm uses "centered" division
\begin{equation}\label{eucl_new}
a = b q + r, \quad -\frac{b}{2} < r \leqslant \frac{b}{2}
\end{equation}
and leads to the following representation of a real number $x$:
\begin{equation}\label{frac_new}
x =\left[a_0; \frac{\varepsilon_1}{a_1},\ldots,
\frac{\varepsilon_l}{a_l}, \ldots \right]= a_0 +
\cfrac{\varepsilon_1}{a_1 + \cfrac{\varepsilon_2}{a_2+ \ldots +
\cfrac{\varepsilon_l}{a_l + \ldots}}}.
\end{equation}
This representation  is known as a continued fraction with minimal
remainders.
Numbers $a_i$ are called partial quotients of fraction (\ref{frac_new}). 
Here $a_0 \in \mathbb{Z},$ $\varepsilon_i = \pm 1$ and
$a_j \geqslant 2$, $a_j + \varepsilon_{j+1} > 2$ for $j \geqslant
1.$ For rational $x$, if $a_s=2$ is the last partial quotient,
then $\varepsilon_s = 1$ for uniqueness of the representation.
Such fractions can be found in the book \cite{Perron} by
O.\,Perron.

Statistical properties of various Euclidean algorithms were
investigated by B.\,Vallee and V.\,Baladi in papers
\cite{vallee1}, \cite{vallee2}, \cite{vallee3}. The most precise
asymptotic formulae for the mean length for the classical
Euclidean algorithm and the centered Euclidean algorithm are
proved in papers \cite{Ustinov}, \cite{Ustinov2} by
A.\,V.\,Ustinov. A similar formula for "by-excess" Euclidean
algorithm was obtained in author's paper \cite{Jabitskaya1}.

For rational  $x$ let us denote by $S(x)$ the sum of partial
quotients of representation (\ref{frac_new}), and put
$$\mathcal{Z}_n: = \left\{x \in \mathbb{Q} \cap [0,1] : S(x) \leqslant n + 1\right\}.$$

In present paper we investigate the limit distribution function
\begin{equation}\label{Z_n}
F(x) := \lim_{n \to \infty} \frac{\sharp \left\{ \xi \in \mathcal{Z}_n: \xi \leqslant x\right\}}
{\sharp \mathcal{Z}_n},\quad x \in [0,1].
\end{equation}
The main result is the following theorem.
\begin{theorem}
Let $x \in [0,1]$, $x$ is represented in the form (\ref{frac_new}), then
%\begin{multline}
%F(x) = a_0 - (-\varepsilon_1)\frac{c}{\lambda^{a_0 + a_1 - 1}} -
%(-\varepsilon_1)(-\varepsilon_2)
%\frac{c}{\lambda^{a_0 + a_1 + a_2 - 1}} - \ldots - \\ -
%\left(\prod_{1 \leqslant j \leqslant l}(-\varepsilon_j)\right)
%\frac{c}{\lambda^{ \sum\limits_{0 \leqslant i \leqslant l}a_i- 1}} - \ldots,
%\end{multline}
\begin{equation}\label{F}
F^{[2]}(x) = a_0 - c \lambda\left(\frac{E_1}{\lambda^{A_1}} +
\frac{E_2}{\lambda^{A_2}} + \ldots +\frac{E_j}{\lambda^{ A_j}} +
\ldots\right),
\end{equation}
where $$E_j = \prod_{1 \leqslant i \leqslant j}(-\varepsilon_i),
\quad A_j = \sum\limits_{0 \leqslant i \leqslant j}a_i, \,\,\, c =
1/(\lambda-1),
$$
and  $\lambda$ is  the unique real root of the equation
$$
\lambda^3 - \lambda^2 - \lambda - 1 = 0.
$$

For rational $x$ the sum in formula (\ref{F}) is finite.
\end{theorem}

In this paper we also prove
\begin{theorem}
Let for $x \in[0,1]$ the derivative $F'(x)$ (finite or infinite)
exists. Then either $F'(x) = 0$ or $F'(x) = \infty$.
\end{theorem}
As function $F(x)$ is monotonic, then by Lebesgue's theorem, the derivative
$F'(x)$ exists and is finite
almost everywhere (in the sense of Lebesgue measure).
That is why $F'(x) = 0$ almost everywhere.
In other words, $F(x)$ is a singular function.

In the proof of Theorem 1 we need the following result.

\begin{Prop}
For $x \in [0, 1/2]$ function $F(x)$ satisfies the following
functional equation
\begin{equation}\label{0}
F(1 - x) = 1 - \frac{F(x)}{\lambda}.
\end{equation}
\end{Prop}

The proof of Proposition 1 is given in section 3. Theorem 1 uses
Proposition 1 and it's proof is given in section 4. The proof of
Theorem 2 is in section 5.
%The rest of the article is in Russian. English version should appear in
%mathematical notes.

\section{Properties of a continued fraction with minimal remainders}

It follows immediately from the definition of continued
fraction with minimal remainders (\ref{frac_new}), that
\begin{itemize}
\item $a_i \geqslant 2$ for $i \geqslant 1$,
\item If $a_i=2$, then $\varepsilon_{i+1} = 1$ for $i \geqslant 1$,
\item If the last partial quotient $a_l = 2$, then $\varepsilon_l = 1$.
\end{itemize}

Let $x = [b_0; b_1, \ldots, b_s, \ldots]$ is represented in the
form of ordinary continued fraction (\ref{ordcf}). We describe the
algorithm for  converting this fraction into a fraction of the
form (\ref{frac_new}) (see   \cite{Perron}).

Fraction (\ref{ordcf}) is constructed by the classical Euclidean
algorithm
$$
r_{0} = \frac{b}{a}, \quad r_{i+1} = \frac{1}{r_{i}} -  b_{i},
\quad b_i = \left[\frac{1}{r_{i}} \right], \quad 0 \leqslant r_{i}
< 1.$$ The remainder $r_{i+1}$ is less then $\cfrac{1}{2}$ if and
only if $b_{i+1} > 1$. So while $b_{i+1} > 1$ partial quotients of
the classical Euclidean algorithm coincide with partial quotients
of "centered" Euclidean algorithm.

For the first $i$ such that $b_{i+1} = 1$, we use the identity
$$
b_{i} + \cfrac{1}{1 + \cfrac{1}{b_{i+2} + \alpha}} = b_{i} + 1 -
\cfrac{1}{b_{i+2} + 1 + \alpha}, \quad \alpha \geqslant 0.
$$
And since $\displaystyle \frac{1}{b_{i+1} + 1} < \frac 1 2$, we have
$$
\left[b_0; b_1, \ldots, b_i, 1 , b_{i+2}\right] =
\left[b_0; \cfrac{1}{b_1},
\ldots, \cfrac{1}{b_i+1}, \cfrac{-1}{b_{i+2}+1} \right]
$$
Then we apply the same procedure to the "tail"
$$
b_{i+2} + 1 + \cfrac{1}{b_{i+3} +
\ldots + \cfrac{1}{b_s + \ldots}}
$$
of the fraction (\ref{ordcf}).

We define the convergents of the continued fraction with minimal
remainders of the number $\displaystyle x = \left[a_0;
\frac{\varepsilon_1}{a_1},\ldots, \frac{\varepsilon_l}{a_l},
\ldots \right]$ as
$$
\frac{P_n}{Q_n} = \left[a_0; \frac{\varepsilon_1}{a_1},\ldots, 
\frac{\varepsilon_n}{a_n}\right],\quad (P_n,Q_n) = 1,
\quad n \geqslant 0.
$$
To get a recurrence formulas for $P_n/Q_n$, $n \geqslant 0$ we put formally
$$\cfrac{P_{-1}}{Q_{-1}} = \frac{1}{0}, \quad
\cfrac{P_{0}}{Q_{0}} = \cfrac{a_0}{1}.$$ Then for $\varepsilon_{n+1} =1$
we have
%If $(n-1)$-th and $n$-th convergents are located by different sides from
%number $\cfrac{a}{b}$ on the number scale, then
$$\cfrac{P_{n+1}}{Q_{n+1}} = \cfrac{a_n P_{n} + P_{n-1}}{a_n Q_{n} +
Q_{n-1}},$$
otherwise
$$\cfrac{P_{n+1}}{Q_{n+1}} =\cfrac{a_n P_{n} - P_{n-1}}
{a_n Q_{n} - Q_{n-1}}.$$

\section{Definition and properties of sets $\mathcal{Z}_n$}

We define a sequence of sets $\mathcal{X}_k$ by
$$\mathcal{X}_k = \left\{x \in \mathbb{Q} \cap [0,1]: S(x) =
k+1\right\}, \quad n \geqslant 1.$$
It is clear that
$$\mathcal{Z}_n = \mathop{\cup}\limits_{1 \leqslant k \leqslant n} \mathcal{X}_k,$$
where $\mathcal{Z}_n$ is defined by (\ref{Z_n}).
Suppose that the elements of $\mathcal{Z}_k$ are
arranged in the increasing order. The number of elements of
$\mathcal{Z}_n$, $\mathcal{X}_n$ we denote by $Z_n$, $X_n$
correspondingly.

Particularly, $\mathcal{X}_1 =\left\{\cfrac 1 2\right\}$,
$\mathcal{X}_2 =\left\{\cfrac 1 3\right\}$, $\mathcal{X}_3
=\left\{\cfrac 1 4, \cfrac 2 5, \cfrac 2 3\right\}$,
$\mathcal{X}_4 =\left\{\cfrac 1 5, \cfrac 2 7, \cfrac 3 7, \cfrac
3 5, \cfrac 3 4\right\}$. So $X_1 = X_2 = 1$, $X_3 = 3$, $X_4 =
5$.

\begin{lemma}
For $n \geqslant 1$ we have
$$X_{n+3} = X_{n+2} + X_{n+1} + X_{n}.$$
\end{lemma}
{\bf Proof.} We construct one-to-one correspondence $\Phi$ between
elements of sets $\mathcal{X}_{n+2} \cup \mathcal{X}_{n+1} \cup
\mathcal{X}_{n}$ and $\mathcal{X}_{n+3}$.

Let $x \in \mathcal{X}_{n+2} \cup \mathcal{X}_{n+1} \cup \mathcal{X}_{n}$,
$x = \left[a_0; \cfrac{\varepsilon_1}{a_1},\ldots, \cfrac{\varepsilon_l}{a_l} \right]$,
we define $\Phi(x):\mathcal{X}_{n+2} \cup \mathcal{X}_{n+1} \cup \mathcal{X}_{n}
\to \mathcal{X}_{n+3}$ in the following way:
\begin{itemize}
\item If $x \in \mathcal{X}_{n+2}$, then
$$\Phi(x)=\left[a_0; \cfrac{\varepsilon_1}{a_1},
\ldots, \cfrac{\varepsilon_{l-1}}{a_{l-1}}, \cfrac{\varepsilon_l}{a_l + 1} \right] \in \mathcal{X}_{n+3}.$$
\item If $x \in \mathcal{X}_{n+1}$, then
$$\Phi(x)=\left[a_0; \cfrac{\varepsilon_1}{a_1},
\ldots, \cfrac{\varepsilon_l}{a_l}, \cfrac{1}{2} \right] \in \mathcal{X}_{n+3}.$$
\item If $x \in \mathcal{X}_{n}$ and $a_{l} > 2$, then
$$\Phi(x)=\left[a_0; \cfrac{\varepsilon_1}{a_1},
\ldots, \cfrac{\varepsilon_l}{a_l}, \cfrac{-1}{3} \right] \in \mathcal{X}_{n+3}.$$
\item If $x \in \mathcal{X}_{n}$ and $a_{l} = 2$, then
$$\Phi(x)=\left[a_0; \cfrac{\varepsilon_1}{a_1}, \ldots, \cfrac{\varepsilon_{l-2}}{a_{l-2}},
\cfrac{\varepsilon_{l-1}}{a_{l-1}+1}, \cfrac{-1}{2}, \cfrac{1}{2} \right]
\in \mathcal{X}_{n+3}.$$
\end{itemize}

The correspondence $\Phi(x)$ is injective by the construction. Let
us show that it is surjective. For any $y \in \mathcal{X}_{n+3}$,
$y=\left[a_0; \cfrac{\varepsilon_1}{a_1}, \ldots,
\cfrac{\varepsilon_l}{a_l} \right]$ we find the preimage $x$ of
$y$.

\begin{itemize}
\item If $a_l > 3$ or $a_l = 3$ and $\varepsilon_l =1$ then
$$x = \left[a_0; \cfrac{\varepsilon_1}{a_1},
\ldots, \cfrac{\varepsilon_{l-1}}{a_{l-1}},\cfrac{\varepsilon_l}{a_{l}-1} \right] \in \mathcal{X}_{n+2}.$$
\item If $a_l = 2$ and either $a_{l-1} > 2$ or $a_{l-1} = 2$ and $\varepsilon_{l-1} =1$,
then $$x = \left[a_0; \cfrac{\varepsilon_1}{a_1},\ldots,
\cfrac{\varepsilon_{l-1}}{a_{l-1}}\right] \in \mathcal{X}_{n+1}.$$
\item If $a_l = 3$, $\varepsilon_l =-1$, then $a_{l-1} > 2$,
therefore $$x = \left[a_0; \cfrac{\varepsilon_1}{a_1},
\ldots, \cfrac{\varepsilon_{l-1}}{a_{l-1}} \right] \in \mathcal{X}_{n}.$$
\item If $a_l = a_{l-1} = 2$, $\varepsilon_{l-1} =-1$, then $a_{l-2} > 2$,
therefore $$x = \left[a_0; \cfrac{\varepsilon_1}{a_1}, \ldots,
\cfrac{\varepsilon_{l-3}}{a_{l-3}},
\cfrac{\varepsilon_{l-2}}{a_{l-2}-1}, \cfrac{1}{2} \right]
\in \mathcal{X}_{n}.$$
\end{itemize}

Lemma is proved.
 \hfill $\blacksquare$

\begin{sled}
For $n \geqslant 1$ we have
\begin{equation}\label{z_n}
Z_{n+3} = Z_{n+2} + Z_{n+1} + Z_{n} + 2.
\end{equation}
\end{sled}
{\bf Proof.} By the definition of $\mathcal{Z}_n$ and Lemma 1, we
get
\begin{multline*}
Z_{n+2} + Z_{n+1} + Z_{n} = \\ =
\left(X_{1} + \ldots + X_{n+2}\right)  + \left(X_{1} + \ldots + X_{n+1}\right) +
\left(X_{1} + \ldots + X_{n} \right)= \\ = X_{1} + X_{2} + X_3 + X_{4} + \ldots +
X_{n+3} + \left(X_{1} - X_{3}\right)=  Z_{n+3} - 2.$$
\end{multline*}
\hfill $\blacksquare$

We remind the definition of the  Stern-Brocot sequences
$\mathcal{F}_n$, $n = 0,1,2,\ldots$. Consider two-point set
$\displaystyle \mathcal{F}_0 = \left\{ \frac 0 1 , \frac 1 1
\right\}$. Let $n \geqslant 0$ and
$$\mathcal{F}_n =
\left\{0 = x_{0,n} < x_{1,n}< \ldots < x_{N(n),n} = 1 \right\},$$
where $x_{j,n} = p_{j,n}/q_{j,n}$, $(p_{j,n},q_{j,n}) = 1$,
$j = 0,\ldots, N(n)$ and $N(n) = 2^n+1$.
Then
$$\mathcal{F}_{n+1} =\mathcal{F}_{n} \cup Q_{n+1}$$
with
$$Q_{n+1} = \left\{x_{j-1,n}\oplus x_{j,n},\quad j = 1, \ldots, N(n)\right\}.$$
Here
$$\frac{a}{b} \oplus \frac{c}{d} = \frac{a+b}{c+d}$$
is the mediant of fractions $\displaystyle\frac{a}{b}$ and
$\displaystyle\frac{c}{d}$.

So the first sequences are  $\displaystyle {Q}_1 =\left\{\frac 1
2\right\}$, $\displaystyle {Q}_2 =\left\{\frac 1 3,
\frac{2}{3}\right\}$, $\displaystyle {Q}_3 =\left\{\frac 1 4,
\frac 2 5, \frac 3 5, \frac{3}{4}\right\}$. It is clear that for
any rational number $q$ there exists such number $n$ that $q \in
{Q}_n$. Note that sum $S^{[0]}(x)$ of partial quotients of
ordinary continued fraction of a number $x \in {Q}_n$ equals to
$n+1$. Formula (\ref{F_n}) gives an equivalent definition of $\mathcal{F}_n$.

It is convenient to represent sequences $\mathcal{F}_n$ by means
of the binary tree $\mathcal{D}^{[0]}$~(Figure~1). This tree is
called Stern-Brocot's tree.
\begin{figure}[h]
\centerline{\includegraphics{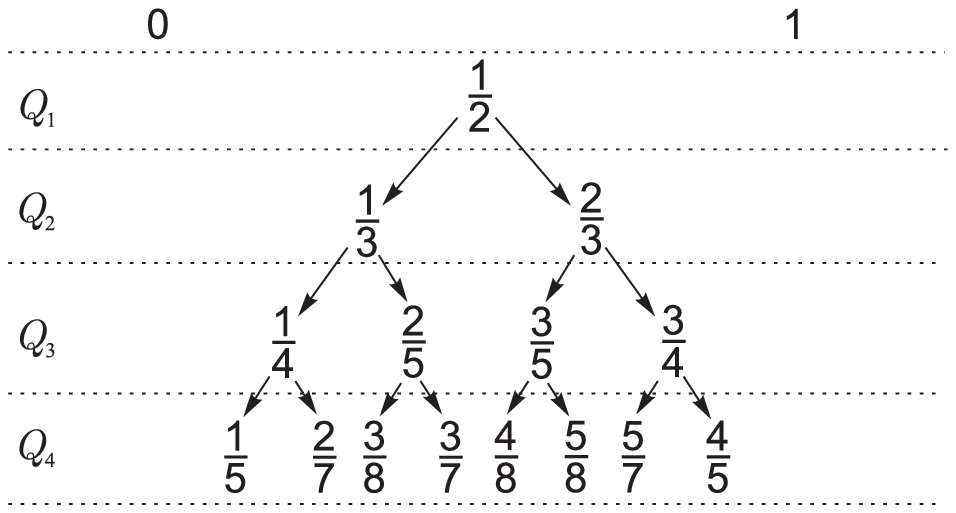}}
\caption{}
\end{figure}
\begin{figure}[h!]
\centerline{\includegraphics{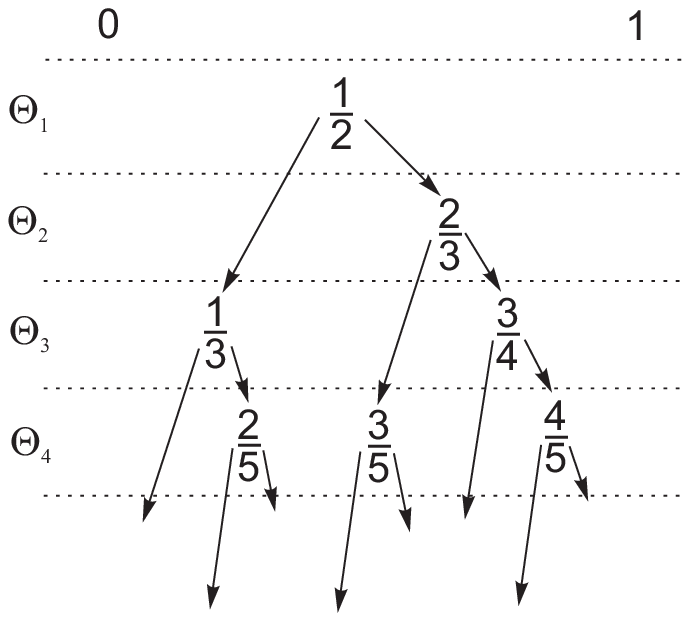}}
\caption{}
\end{figure}
Nodes of the tree are labeled by rationals from $(0,1)$ and partitioned into
levels by the following rule:
%By $\mathcal{D}^{[0]}$ (Figure~1) we denote tree with nodes distributed by levels
$n$-th level consists of nodes labeled by numbers $x$, such that
$S^{[0]}(x) = n+1$ (i.e. $n$-th level consists of nodes, labeled by numbers
from ${Q}_n$).

It is possible to distribute nodes of the tree into levels by
another way. For example, we can use such a rule: $n$-th level
consists of nodes labeled by numbers $x$, such that sum
$S^{[1]}(x)$ of partial quotients of regular reduced continued
fraction of number $x$ equals $n+1$. Then we get tree
$\mathcal{D}^{[1]}$~(Figure~2) from paper \cite{jabitskaya2}.

Now let us distribute nodes of the tree into levels by the
following rule:
%In present paper we need tree $\mathcal{D}$, which nodes are distributed
%by levels in the following way:
$n$-th level consists of nodes labeled by numbers $x$, such that
$S(x) = n+1$ (i.e. $x \in \mathcal{X}_n$).
\begin{figure}[h]
\centerline{\includegraphics{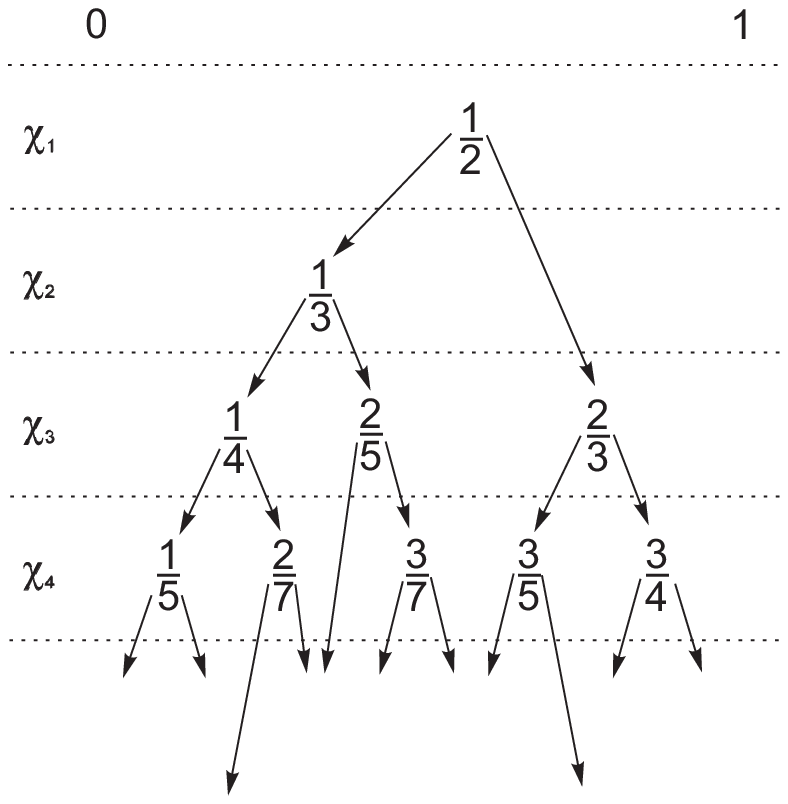}}
\caption{}
\end{figure}
We denote this tree by $\mathcal{D}$~(Figure~3).

Any node $\xi$ of the tree $\mathcal{D}$ is a root of subtree,
which we denote by $\mathcal{D}^{(\xi)}$. Nodes of
$\mathcal{D}^{(\xi)}$ are also partitioned into levels: $\xi$
itself belongs to level 1, and a node of the tree
$\mathcal{D}^{(\xi)}$, labeled by number $x$ belongs to the level
$S(x) - S(\xi) + 1$ in the tree $\mathcal{D}^{(\xi)}$. The number
of nodes of $\mathcal{D}^{(\xi)}$ from level 1 to level $n$   we
denote by $D^{(\xi)}_n$.

Let us consider more detailed structure of the tree $\mathcal{D}$.
From every node $\xi$ of $\mathcal{D}$ we issue two arrows: the
left one and the right one. The left one goes to the node labeled
by $x^l$ and the right one goes to node labeled by $x^r$. Note
that if $\xi = x \oplus y$, where $x$, $y$ are consecutive
elements of $\mathcal{F}_n$, then $\xi^l = x \oplus \xi$, $\xi^r =
\xi \oplus y$.
\begin{zam}\label{zam}
Let $\xi \in \mathcal{X}_n$ and $\xi = \left[a_0; \cfrac{\varepsilon_1}{a_1},
\ldots, \cfrac{\varepsilon_l}{a_{l}} \right] = x \oplus y$, where $x, y$~---
neighboring elements of $\mathcal{F}_n$, $S(x) < S(y)$.
If $a_l >2$, then
$$x \oplus \xi = \left[a_0; \cfrac{\varepsilon_1}{a_1},
\ldots, \cfrac{\varepsilon_{l-1}}{a_{l-1}}, \cfrac{\varepsilon_l}{a_{l}+1} \right] \in X_{n+1},$$
$$y \oplus \xi = \left[a_0; \cfrac{\varepsilon_1}{a_1},\ldots,\cfrac{\varepsilon_{l-1}}{a_{l-1}}, 
\cfrac{\varepsilon_l}{a_{l}-1},
\cfrac{1}{2} \right]\in X_{n+1},$$
%i.\,e. $\xi^l, \xi^r \in X_{n+1}$.
If $a_l = 2$, then $\varepsilon_l = 1$ and
$$x \oplus \xi = \left[a_0; \cfrac{\varepsilon_1}{a_1},
\ldots, \cfrac{\varepsilon_{l-1}}{a_{l-1}}, \cfrac{1}{3} \right] \in X_{n+1},$$
$$y \oplus \xi = \left[a_0; \cfrac{\varepsilon_1}{a_1},\ldots,
\cfrac{\varepsilon_{l-1}}{a_{l-1}},\cfrac{\varepsilon_{l-1}}{a_{l-1}+1},
\cfrac{-1}{3} \right] \in X_{n+2}.$$
\end{zam}
From Remark \ref{zam} we deduce the following statement.
\begin{lemma}\label{lemma1}
Let $\xi = \left[a_0; \cfrac{\varepsilon_1}{a_1},
\ldots, \cfrac{\varepsilon_l}{a_{l}} \right]$, then
\begin{equation}
D^{(\xi)}_n =
\begin{cases}
D^{(1/2)}_n, \quad \text{if } a_l = 2\\
D^{(1/3)}_n, \quad \text{if } a_l > 2.
\end{cases}
\end{equation}
\end{lemma}

Note that $D^{(1/2)}_n = Z_n$. For  brevity we put $Y_n
=D^{(1/3)}_n$. So $Y_1 = 1$, $Y_2 = 2$, $Y_3 = 3$.

For $Z_n$ we have recurrence formula (\ref{z_n}). It is easy to
prove a similar formula for $Y_n$:
\begin{equation}\label{y_n}
Y_{n+3} = Y_{n+2} + Y_{n+1} + Y_{n} + 2.
\end{equation}

\begin{lemma}

Let $\lambda$ be the unique real root of the equation
\begin{equation}\label{haract}
\lambda^3 - \lambda^2 - \lambda - 1 = 0,
\end{equation}
and $c = 1/(\lambda-1)$. Then

$$\lim_{n \to \infty} \frac{Y_n}{Y_{n+1}} =
\lim_{n \to \infty} \frac{Z_n}{Z_{n+1}} = \lambda, \quad \lim_{n
\to \infty} \frac{Y_n}{Z_n} = c.$$

\end{lemma}
%Let us find explicit formula for $Y_n$, $Z_n$.
{\bf Proof.} Equation (\ref{y_n}) can be reduced to a homogeneous
by the substitution  $Y'_n = Y_n +1$:
$$Y'_{n+3} = Y'_{n+2} + Y'_{n+1} + Y'_{n}.$$
The characteristic equation is
\begin{equation*}
\lambda^3 - \lambda^2 - \lambda = 1.
\end{equation*}
This equation has the unique real root $\lambda \approx 1,839292$
and two complex roots $\lambda_1$, $\lambda_2$, such that
$|\lambda_1| = |\lambda_2| < 1$. So
\begin{gather*}Y_n + 1 = Y'_n = C_1 \lambda^n + C_2 \lambda^n_2 +  C_3 \lambda^n_3,\\
Z_n + 1 =Z'_n = D_1 \lambda^n + D_2 \lambda^n_2 +  D_3
\lambda^n_3,
\end{gather*}
with certain constants $C_1$, $C_2$, $C_3$, $D_1$, $D_2$, $D_3$.
%The ratio $C_1/D_1$ we denote by $c$.
Put $c = C_1/D_1$.
%We also note that
From construction of the tree $\mathcal{D}$ it is clear that
$$Z_n = Y_{n-1} + Y_{n-2} +1.$$
Dividing both parts of this equality by $Z_n$ and taking the limit
we get
\begin{multline*}
1 = \lim_{n \to \infty} \frac{C_1 \left(\lambda^{n-1}_1 + \lambda^{n-2}_1\right) +
C_2 \left(\lambda^{n-1}_2 + \lambda^{n-2}_2\right) +
C_3 \left(\lambda^{n-1}_3 + \lambda^{n-2}_3\right) + 1}
{D_1 \lambda^n_1 + D_2 \lambda^n_2 +  D_3 \lambda^n_3 + 1} = \\
= \frac{C_1 + C_1\lambda}{D_1 \lambda^2}.
\end{multline*}
Since $\lambda$ is the root of equation (\ref{haract}), we get the
following relation between  $c$ and $\lambda$:
\begin{equation}\label{c}
c = \frac{\lambda^2}{1 + \lambda} = \frac{1}{\lambda - 1}\approx 1,1915.
\end{equation}

Lemma is proved.
 \hfill $\blacksquare$

\section{Properties of the limit distribution function $F(x)$ of
sequence $\mathcal{Z}_n$}
In this section we prove some auxiliary results about function $F(x)$.

\begin{lemma}\label{lemma2}
Let $x$, $y$ be consecutive elements of the sequence
$\mathcal{Z}_n$. Then
$$
F(y) - F(x) = \lim_{n \to \infty} \frac{D^{(x \oplus y)}_{n+2-S(x \oplus y)}
}{D^{(1/2)}_n}.
$$
\end{lemma}
{\bf Proof.} Note that nodes of tree $\mathcal{D}^{(x \oplus y)}$
are labeled exactly by the numbers from the set $\left\{ \xi \in
\mathbb{Q}: x < \xi < y\right\}.$
%It is enough to note that the number of nodes of the tree
%$\mathcal{D}^{(x \oplus y)}$
%labeled by numbers from $\mathcal{Z}_n$ is exactly
%$D^{(x \oplus y)}_{n+2-S(x \oplus y)}$.
So
$$
F(y) - F(x) =\lim_{n \to \infty} \frac{\sharp \left\{ \xi \in \mathcal{Z}_n:
x < \xi \leqslant y\right\}}{{Z}_n} =
\lim_{n \to \infty} \frac{D^{(x \oplus y)}_{n+2-S(x \oplus y)}}{D^{(1/2)}_n}
.$$
\hfill $\blacksquare$
\begin{lemma}\label{ratio}
Let $x$, $y$ be consecutive elements of $Z_n$, $S(x) < S(y)$, and
let $a_l$ be the last partial quotient in continued fraction with
minimal remainders
representation %(\ref{frac_new})
of the number $x \oplus y$.

If $a_l = 2$, then
\begin{gather}\label{ration1}
F(x \oplus y) - F(x) = \frac{c}{\lambda} \left(F(y) - F(x)\right),\\
\label{ration1.1} F(x \oplus y) - F(y) = \frac{c}{\lambda^2}
\left(F(x) - F(y)\right).
\end{gather}
If $a_l > 2$, then
\begin{gather}\label{ration2}
F(x \oplus y) - F(x) = \frac{1}{\lambda} \left(F(y) - F(x)\right),\\
\label{ration2.2}
F(x \oplus y) - F(y) = \frac{1}{c\lambda} \left(F(x) - F(y)\right).
\end{gather}
\end{lemma}
{\bf Proof.} We suppose that $x < y$ (in case $x > y$ the proof is
similar). According to Lemma \ref{lemma2} one has
$$
F(x \oplus y) - F(x) = \lim_{n \to \infty}
\frac{D^{((x \oplus y)^l)}_{n+2 - S((x \oplus y)^l)}}{D_n}
$$
$$
F(y) - F(x\oplus y) = \lim_{n \to \infty}
\frac{D^{((x \oplus y)^r)}_{n+2 - S((x \oplus y)^r)}}{D_n}
$$

By Remark \ref{zam}, if $a_l = 2$, then the last partial quotients
of continued fractions with minimal remainders of numbers $(x
\oplus y)^l$, $(x \oplus y)^r$ are greater then 2 and
%������� $(x \oplus y)^l$ � $\mathcal{D}$ �� ������� ������ ������ $(x \oplus y)^r$.
$$
S((x \oplus y)^l) = S(x \oplus y) + 1, \quad 
S((x \oplus y)^r) = S(x \oplus y) + 2.
$$
So
\begin{equation}
\frac{F(x \oplus y) - F(x)}{F(y) - F(x\oplus y)} =
\lim_{n \to \infty} \frac{D^{((x \oplus y)^l)}_{n+2 -S((x \oplus y)^l)}}
{D^{((x \oplus y)^r)}_{n+2 -S((x \oplus y)^r)}} =
\lim_{n \to \infty} \frac{Y_{n+1-S(x \oplus y)}}{Y_{n-S(x \oplus y)}} =
{\lambda}.
\end{equation}
i.e. $F(x \oplus y)$ divides the segment $\left[F(x), F(y)\right]$
in the relation $\lambda : 1$. Taking into account that
$$\cfrac{c}{\lambda^2} + \cfrac{c}{\lambda} =1$$
we get formulas (\ref{ration1}),  (\ref{ration1.1}).

If $a_l > 2$, then by Remark \ref{zam} the last partial quotient
of the continued fraction with minimal remainders for the  number
$(x \oplus y)^r$ is 2. For the number $(x \oplus y)^l$ the last
partial quotient  is greater then 2 and
%������ $(x \oplus y)^l$ � $(x \oplus y)^r$ � $D$ ����� ���������.
$$
S((x \oplus y)^l) = S((x \oplus y)^r) = S(x \oplus y) +1.
$$
That is why
\begin{equation}
\frac{F(x \oplus y) - F(x)}{F(y) - F(x\oplus y)} =
\lim_{n \to \infty} \frac{D^{((x \oplus y)^l)}_{n+2 - S(x \oplus y)^l)}}
{D^{((x \oplus y)^r)}_{n+2 - S(x \oplus y)^r)}} =
\lim_{n \to \infty} \frac{Y_{n+1 -S(x \oplus y)} }{Z_{n+1 -S(x \oplus y)}} =
{c}.
\end{equation}
i.e. the number $F(x \oplus y)$ divides the segment $\left[F(x),
F(y)\right]$ in the relation $c:1$. Taking into account that
$$\cfrac{1}{\lambda} + \cfrac{1}{c\lambda} =1$$
we get formulas (\ref{ration2}), (\ref{ration2.2}). Lemma is
proved. $\blacksquare$

Now we are able to prove Proposition 1.  Suppose $x \in[0, 1/2]$.
By the definition
$$
F(x) = \lim_{n \to \infty} \frac{\sharp \left\{ \xi \in \mathcal{Z}_n:
\xi \leqslant x\right\}}{{Z}_n}.
$$
So
$$
\lim_{n \to \infty} \frac{\sharp \left\{ \xi \in \mathcal{Z}_{n-1}: \xi 
\leqslant x\right\}}{{Z}_{n-1}} \frac{Z_{n-1}}{Z_{n}} = \frac{F^{[2]}(x)}{\lambda}.
$$
Taking into account that $S(1-\xi) = 1 + S(\xi)$ for $\xi \in \mathbb{Q}
\cap [0, 1/2)$, for $x \in [0, 1/2]$ we have 
\begin{multline*}
\left\{ \xi \in \mathcal{Z}_{n-1}: \xi < x\right\} =
\left\{ 1-\xi = \eta \in \mathcal{Z}_{n}: 1-\eta < x\right\} = \\=
\left\{ \eta \in \mathcal{Z}_{n}: \eta > 1-x\right\}.
\end{multline*}
So we get
\begin{multline*}
\frac{F^{[2]}(x)}{\lambda} + {F^{[2]}(1-x)} = \\=
\lim_{n \to \infty} \frac{\sharp \left\{ \xi \in \mathcal{Z}_{n-1}: 
\xi \leqslant x\right\}}{{Z}_n} +
\lim_{n \to \infty} \frac{\sharp \left\{ \xi \in \mathcal{Z}_{n}: \xi 
\leqslant 1-x\right\}}{Z_n} = \\ =
\lim_{n \to \infty} \frac{\sharp \left\{ \eta \in \mathcal{Z}_n: 
\eta > 1-x\right\} \cup \left\{ \xi \in \mathcal{Z}_{n}: 
\xi \leqslant 1-x\right\}}{{Z}_n} =1.
\end{multline*}
This equality proofs formula (\ref{0}).
\hfill $\blacksquare$
\section{Proof of Theorem 1}
Let us prove the theorem for rational $x \in [0,1/2]$ by induction
on $S(x)$. The equality
$$F(1/a_1) = \frac{c}{\lambda^{a_1-1}},$$
follows from formula (\ref{ration1}), since
$1/a_1 = \underbrace{0 \oplus \ldots \oplus 0 \oplus }_{(a_1-1)\,\, \text{times}} 1$.

Suppose that the formula
(\ref{F}) is proved for
$$x = \left[0;\frac{1}{a_1},\frac{\varepsilon_2}{a_2},
\ldots,\frac{\varepsilon_m}{a_m}\right].$$
Then it is enough to prove it for
$$y = \left[0;\frac{1}{a_1},\frac{\varepsilon_2}{a_2}, \ldots,
\frac{\varepsilon_{m-1}}{a_{m-1}},\frac{\varepsilon_m}{a_m+1}\right]$$
and for
$$z = \left[0;\frac{1}{a_1},\frac{\varepsilon_2}{a_2}, \ldots,
\frac{\varepsilon_m}{a_m-1}, \frac{1}{2}\right],\quad \text{if}\quad a_m>2,$$
$$w = \left[0;\frac{1}{a_1},\frac{\varepsilon_2}{a_2}, \ldots,
\frac{\varepsilon_{m-2}}{a_{m-2}},\frac{\varepsilon_{m-1}}{a_{m-1} + 1},
\frac{-1}{3}\right],\quad \text{if}\quad a_m = 2.$$

We see that
$$
y = \left[0;\frac{1}{a_1},\frac{\varepsilon_2}{a_2}, \ldots, \frac{\varepsilon_{m-1}}{a_{m-1}}\right] \oplus
\left[0;\frac{1}{a_1},\frac{\varepsilon_2}{a_2}, \ldots, \frac{\varepsilon_{m-1}}{a_{m-1}}, \frac{\varepsilon_{m}}{a_{m}}\right]
$$
and the last partial quotient $a_m+1$ of continued fraction with
minimal remainders expression of number $y$ is greater then $2$.
From (\ref{ration2}) and the inductive assumption  we get
\begin{multline*}
F(y) =
F \left(\left[0;\frac{1}{a_1},\frac{\varepsilon_2}{a_2}, \ldots, \frac{\varepsilon_{m-1}}{a_{m-1}}\right]\right) + \\ +
\frac{1}{\lambda} \left(F\left(\left[0;\frac{1}{a_1},\frac{\varepsilon_2}{a_2}, \ldots, \frac{\varepsilon_{m}}{a_{m}}\right]\right)-
F \left(\left[0;\frac{1}{a_1},\frac{\varepsilon_2}{a_2}, \ldots, \frac{\varepsilon_{m-1}}{a_{m-1}}\right]\right)\right)= \\ =
F \left(\left[0;\frac{1}{a_1},\frac{\varepsilon_2}{a_2}, \ldots, \frac{\varepsilon_{m-1}}{a_{m-1}}\right]\right) -
\frac{1}{\lambda}
%\left(\prod_{1 \leqslant j \leqslant m} (-\varepsilon_j)\right)
c\lambda\frac{E_m}{\lambda^{A_m}}
%\sum\limits_{0 \leqslant i \leqslant m}a_i- 1}}
= \\ =
F \left(\left[0;\frac{1}{a_1},\frac{\varepsilon_2}{a_2}, \ldots,
\frac{\varepsilon_{m-1}}{a_{m-1}}\right]\right) -
%\left(\prod_{1 \leqslant j \leqslant m} (-\varepsilon_j)\right)
c\lambda \frac{E_m}{
\lambda^{%\sum\limits_{0 \leqslant i \leqslant m-1}a_i
A_{m-1}+ (a_m + 1)}}.
\end{multline*}
If $a_m >2$, we must prove the formula (\ref{F}) for $z$. We see
that
$$
z = \left[0;\frac{1}{a_1},\frac{\varepsilon_2}{a_2}, \ldots, \frac{\varepsilon_m}{a_m}\right]
\oplus \left[0;\frac{1}{a_1},\frac{\varepsilon_2}{a_2}, \ldots,
\frac{\varepsilon_{m-1}}{a_{m-1}},\frac{\varepsilon_m}{a_m-1}\right]
$$
and the last partial quotient of number $z$ is 2. So  by
(\ref{ration1}) we have
\begin{multline*}
F(z) =
F \left(\left[0;\frac{1}{a_1}, \frac{\varepsilon_2}{a_2}, \ldots,
\frac{\varepsilon_{m-1}}{a_{m-1}}, \frac{\varepsilon_m}{a_m-1}\right]\right) - \\ -
\frac{c}{\lambda}\left(F\left(\left[0;\frac{1}{a_1},\frac{\varepsilon_2}{a_2}, \ldots,
\frac{\varepsilon_{m-1}}{a_{m-1}}, \frac{\varepsilon_m}{a_m-1}\right]\right) -
F\left(\left[0;\frac{1}{a_1},\frac{\varepsilon_2}{a_2}, \ldots, \frac{\varepsilon_m}{a_m}\right]\right)\right)= \\ =
F \left(\left[0;\frac{1}{a_1},\frac{\varepsilon_2}{a_2}, \ldots,
\frac{\varepsilon_{m-1}}{a_{m-1}}, \frac{1}{a_m}\right]\right) +
%\left(\prod_{1 \leqslant j \leqslant m} (-\varepsilon_j )\right)
\frac{c}{\lambda}c \lambda
\left(\frac{E_m}{\lambda^{A_{m-1} + (a_m - 1)}} -
\frac{E_m}{\lambda^{A_m}}\right) = \\ =
F \left(\left[0;\frac{1}{a_1},\frac{\varepsilon_2}{a_2}, \ldots,\frac{\varepsilon_{m-1}}{a_{m-1}},
\frac{1}{a_m}\right]\right) -
%\left(\prod_{1 \leqslant j \leqslant m + 1} (-\varepsilon_j )\right)
c\lambda
\frac{E_{m}(-1)}{\lambda^{ A_{m-1}
%\sum\limits_{0 \leqslant i \leqslant m-1}a_i
+ (a_m - 1) +2 }}.
\end{multline*}
If $a_m =2$, we must prove formula (\ref{F}) for $w$. We see that
$$
w = \left[0;\frac{1}{a_1},\frac{\varepsilon_2}{a_2}, \ldots,
\frac{\varepsilon_{m-1}}{a_{m-1}},\frac{\varepsilon_{m}}{a_{m}+1}\right] \oplus
\left[0;\frac{1}{a_1},\frac{\varepsilon_2}{a_2}, \ldots,
\frac{\varepsilon_{m}}{a_{m}}, \frac{1}{2}\right]
$$
and the last partial quotient of number $w$ is greater then 2. So
by  (\ref{ration2}) we have
\begin{multline}\label{F(w)}
F^{[2]}(w) = 
F^{[2]} \left(\left[0;\frac{1}{a_1},\frac{\varepsilon_2}{a_2}, \ldots,
\frac{\varepsilon_{m-2}}{a_{m-2}}, \frac{\varepsilon_{m-1}}{a_{m-1}+1}\right]\right) + \\ + 
\frac{1}{\lambda}\left(F^{[2]} \left(\left[0;\frac{1}{a_1},\frac{\varepsilon_2}{a_2}, \ldots,\frac{\varepsilon_{m-1}}{a_{m-1}}, \frac{1}{2}\right]\right) - \right.\\ \left.-
F^{[2]}\left(\left[0;\frac{1}{a_1},\frac{\varepsilon_2}{a_2}, \ldots, 
\frac{\varepsilon_{m-2}}{a_{m-2}},\frac{\varepsilon_{m-1}}{a_{m-1} + 1}\right]
\right)\right),
\end{multline}
As
\begin{multline*}
\left[0;\frac{1}{a_1},\frac{\varepsilon_2}{a_2}, \ldots,
\frac{\varepsilon_{m-1}}{a_{m-1}}, \frac{1}{2}\right]=\\=
\left[0;\frac{1}{a_1},\frac{\varepsilon_2}{a_2}, \ldots, 
\frac{\varepsilon_{m-1}}{a_{m-1}+1}\right]
\oplus \left[0;\frac{1}{a_1},\frac{\varepsilon_2}{a_2}, 
\ldots,\frac{\varepsilon_{m-2}}{a_{m-2}}, 
\frac{\varepsilon_{m-1}}{a_{m-1}}\right],
\end{multline*}
 by  (\ref{ration1.1}) we have
\begin{multline*}
F^{[2]} \left(\left[0;\frac{1}{a_1},\frac{\varepsilon_2}{a_2}, \ldots,\frac{\varepsilon_{m-1}}{a_{m-1}}, \frac{1}{2}\right]\right) - \\-
F^{[2]}\left(\left[0;\frac{1}{a_1},\frac{\varepsilon_2}{a_2}, \ldots, 
\frac{\varepsilon_{m-2}}{a_{m-2}},\frac{\varepsilon_{m-1}}{a_{m-1} + 1}\right]\right) = 
\\ =
\frac{c}{\lambda^2}\left(F^{[2]} \left(\left[0;\frac{1}{a_1},\frac{\varepsilon_2}{a_2}, \ldots, \frac{\varepsilon_{m-1}}{a_{m-1}}\right]\right) \right.- 
\allowdisplaybreaks\\ \left.-
F^{[2]}\left(\left[0;\frac{1}{a_1},\frac{\varepsilon_2}{a_2}, \ldots, 
\frac{\varepsilon_{m-2}}{a_{m-2}},\frac{\varepsilon_{m-1}}{a_{m-1} + 1}\right]\right)\right) = \\ =
-\frac{c}{\lambda^2} c \lambda E_{m-1}
\left(\frac{1}{\lambda^{ \sum\limits_{0 \leqslant i \leqslant m-2}a_i + 
a_{m - 1}}} -  \frac{1}{\lambda^{ \sum\limits_{0 \leqslant i \leqslant m-2}a_i + 
(a_{m - 1} + 1)}}\right)=\\=
-c \lambda E_{m-1}\frac{1}{\lambda^{ \sum\limits_{0 \leqslant i \leqslant m-1}a_i
+3}}.
\end{multline*}
Substituting this result in (\ref{F(w)}), we finally get
\begin{equation*}
F^{[2]}(w) = 
F^{[2]} \left(\left[0;\frac{1}{a_1},\frac{\varepsilon_2}{a_2}, \ldots, 
\frac{\varepsilon_{m-2}}{a_{m-2}},
\frac{\varepsilon_{m-1}}{a_{m-1}+1}\right]\right) -  
c \lambda E_{m-1}\frac{1}{\lambda^{ A_{m-1}+4}}.
\end{equation*}

So we have proven Theorem 1 for rational $x \in [0, 1/2]$. For
rational $x \in (1/2,1]$ it follows from formula (\ref{0}). For
irrational $x\in [0,1]$ we should take into account the
continuity of $F(x)$. $\blacksquare$

\section{Singularity of the function $F(x)$}
%Function $F(x)$ has the same property that functions from the family of
%Tichy-Uitz $g_{\lambda}$.
%\begin{theorem}
%Let for $x \in[0,1]$ the derivative $F'(x)$ (finite or infinite)
%exists. Then either $F'(x) = 0$ or $F'(x) = \infty$.
%\end{theorem}
%Function $F(x)$ is monotonic, so, by Lebesgue's theorem, the derivative
%$F'(x)$ exists and is finite
%almost everywhere (in the sense of Lebesgue measure).
%That is why $F'(x) = 0$ almost everywhere.
%In other words, $F(x)$ is a singular function.
In this section we prove Theorem 2.
First of all let us consider the  case $x \in \mathbb{Q}$.
\begin{lemma}
For rational $x \in [0,1]$ we have
$F'(x) = 0$.
\end{lemma}
Let us prove that the right derivative $F'_{+}(x)$  exists and $F'_{+}(x) = 0$
(for $F'_{-}(x)$ the prove is similar).

Let $x = a/b$, $a, b \in \mathbb{N}$, then there exists such $n$, that
$a/b \in \mathcal{X}_n$.
Denote by $a'/b'$ the right neighbouring  to  $a/b$ element in $\mathcal{Z}_n$.
Sequence of mediants $y_k =\left\{\frac{ka+a'}{kb+b'}\right\}$,
converges to $a/b$ from the right as $k\to \infty$. So for $\xi >x $ sufficiently close to $x$  there exists such $m$ that $x < y_{m+1} \leqslant \xi \leqslant y_{m}$
and so
$$
\frac{\left|F(\xi) - F(a/b)\right|}{\xi - a/b} \leqslant
\frac{F(y_{m}) - F(a/b)}{y_{m+1} - a/b}.
$$
Remind that  $a/b \in \mathcal{Z}_n$, but $a/b \notin \mathcal{Z}_{n-1}$. So
%$a/b = y \oplus a'/b'$ for some $y \in \mathcal{Z}_{n-1}$,then
$S(a/b) > S(a'/b')$. By Lemma \ref{ratio} we see that
\begin{multline*}
F(a/b \oplus a'/b') - F(a/b) \leqslant \max\left(\frac{c}{\lambda^2},
\frac{1}{c\lambda}\right) \left(F(a'/b') - F(a/b)\right) = \\=
\frac{1}{c\lambda}\left(F(a'/b') - F(a/b)\right).
\end{multline*}
Similarly $\displaystyle S\left(\frac{ka+a'}{kb+b'}\right) > S(a/b)$, $k = 1, \ldots, m-1,$ and
\begin{multline*}
F\left(\frac{(k+1)a+a'}{(k+1)b+b'}\right) - F(a/b) \leqslant
\max\left(\frac{c}{\lambda},
\frac{1}{\lambda}\right) \left(F\left(\frac{ka+a'}{kb+b'}\right)
- F(a/b)\right)\leqslant \\\leqslant
\frac{c^k}{\lambda^k}\left(F(a/b \oplus a'/b') - F(a/b)\right).
\end{multline*}
So
\begin{multline}
0\leqslant F'_{+}(x)=\lim_{\xi \to x_+} \frac{F(\xi) - F(x)}{\xi - x}
\leqslant \lim_{m \to \infty} \frac{F(y_m) - F(x)
}{y_{m+1}-x} =\\ = \lim_{m \to \infty} \frac{\displaystyle\frac{c^{m-1}}{\lambda^{m+1}}
(F(a'/b') - F(a/b))}
{\displaystyle\frac{1}{((m+1)b+b')b}} = 0.
\end{multline}
\hfill $\blacksquare$

Let $x \notin \mathbb{Q}$.
Suppose that $F' (x) = a$, where $a$ is finite and $a\not = 0$.
We should prove that it is not possible.
We shall use  the Stern-Brocot sequences $\mathcal{F}_n$.

Given $n$ 
we can find two consecutive elements $p_n/q_n<p'_n/q'_n$ from the set $\mathcal{F}_n$
such that $p_n/q_n<x<p'_n/q'_n$.  In such way we obtain an
infinite sequence of pairs of elements $\{p_n/q_n, p'_n/q'_n\}$,
converging to $x$ from the left and from the right correspondingly.
So
\begin{equation*}
\lim_{n \to \infty} \frac{F (p'_n/q'_n) - F (p_n/q_n)}{p'_n/q'_n-p_n/q_n}
= a,
\end{equation*}
Put
\begin{equation*}
G_n(x) =\frac{F (p'_{n+1}/q'_{n+1}) - F (p_{n+1}/q_{n+1})}{F (p'_{n}/q'_{n}) -
F (p_{n}/q_{n})}
\frac{q_{n+1}q'_{n+1}}{q_n q'_n}.
\end{equation*}
Then
\begin{equation}\label{G_n}
G_n(x) =
\frac{F (p'_{n+1}/q'_{n+1}) - F (p_{n+1}/q_{n+1})}{p'_{n+1}/q'_{n+1}-p_{n+1}/q_{n+1}}
\frac{p'_n/q'_n-p_n/q_n}{F (p'_{n}/q'_{n}) - F (p_{n}/q_{n})}
\xrightarrow[n \to \infty]{} 1.
\end{equation}

It is clear that 
 if $x \in
(p_n/q_n, p_n/q_n \oplus p'_n/q'_n )$ then
the pair $\{p_{n+1}/q_{n+1},p'_{n+1}/q'_{n+1}\}$
coincides with $\{p_n/q_n, p_n/q_n \oplus p'_n/q'_n \}$.
Also   if
$x \in (p_n/q_n \oplus p'_n/q'_n, p'_n/q'_n)$ then
the pair $\{p_{n+1}/q_{n+1},p'_{n+1}/q'_{n+1}\}$
coincides with
$\{p_n/q_n \oplus p'_n/q'_n, p'_n/q'_n\}$.

Note that $p_n/q_n$, $p'_n/q'_n$  are always among the  intermediate
and convergent fractions to $x$ in the sense of ordinary continued fraction.

Let us show that for an irrational $x$ one  can find an infinite subsequence
$\{p_{n_k}/q_{n_k}, p'_{n_k}/q'_{n_k}\}$ of the sequence
$\{p_n/q_n, p'_n/q'_n\}$ with the following property: 
%for any $n$ 
the last partial quotient in the
  continued fraction with minimal remainders
expression of $p_{n_k}/q_{n_k}\oplus p'_{n_k}/q'_{n_k}$
 is equal to 2.

Let 
$$\displaystyle x = \left[0;\frac{1}{a_1},\frac{\varepsilon_2}{a_2}, \ldots,
\frac{\varepsilon_{m}}{a_{m}}, \ldots\right].
$$
Then either
\begin{equation*}\label{e1}
p_n/q_n \oplus p'_n/q'_n = \left[0;\frac{1}{a_1},\frac{\varepsilon_2}{a_2},
\ldots, \frac{\varepsilon_{m-1}}{a_{m-1}},\frac{1}{2}\right],
\end{equation*}
or
\begin{equation*}\label{e2}
p_n/q_n \oplus p'_n/q'_n = \left[0;\frac{1}{a_1},\frac{\varepsilon_2}{a_2},
\ldots, \frac{\varepsilon_{m-1}}{a_{m-1}},\frac{\varepsilon_{m}}{b_{m}}\right],
\end{equation*}
for some natural $m$, $b_m \leqslant a_m +1$.
In the first case the pair $\{p_{n}/q_{n}, p'_{n}/q'_{n}\}$
satisfies the necessary property.

In the other case
we consider the pair $\{p_{n+1}/q_{n+1},
p'_{n+1}/q'_{n+1}\}$ and for $p_{n+1}/q_{n+1}\oplus p'_{n+1}/q'_{n+1} \in
\{p_{n+2}/q_{n+2},p'_{n+2}/q'_{n+2}\}$ by remark \ref{zam} we have 
two possibilities: either
$$
 \left[0;\frac{1}{a_1},
\frac{\varepsilon_2}{a_2}, \ldots, \frac{\varepsilon_{m-1}}{a_{m-1}},
\frac{\varepsilon_{m}}{b_{m}-1}, \cfrac{1}{2}\right]
$$
 or
$$ \left[0;\frac{1}{a_1},\frac{\varepsilon_2}{a_2},
\ldots, \frac{\varepsilon_{m-1}}{a_{m-1}},\frac{\varepsilon_{m}}{b_{m}+1}\right] $$
where 
$b_m + 1 \leqslant a_m +1$. But the last inequality can occur only finitely many times. 
So for any $m$ we can find a pair $\{p_{n_k}/q_{n_k}, p'_{n_k}/q'_{n_k}\}$
satisfying the necessary property.

In the sequel we consider such a subsequence 
$\{p_{n_k}/q_{n_k}, p'_{n_k}/q'_{n_k}\}$.
From  (\ref{G_n}), we see that
\begin{equation}\label{G_n_k}
\lim_{k \to \infty}G_{n_k}(x) = 1.
\end{equation}
For a fixed $k$ we consider following cases:

\textbf{Case 1}. $x \in (p_{n_k}/q_{n_k}, p_{n_k}/q_{n_k} \oplus p'_{n_k}/q'_{n_k} )$

and 

\textbf{Case 2}. $x \in (p_{n_k}/q_{n_k} \oplus p'_{n_k}/q'_{n_k}, p'_{n_k}/q'_{n_k})$.

Consider \textbf{Case 1}. 
As the last partial quotient of $p_n/q_n \oplus p'_n/q'_n$
is equal to $2$, Lemma \ref{ratio} leads to
\begin{multline}\label{case_1}
\frac{F (p'_{n_k+1}/q'_{n_k+1}) - F (p_{n_k+1}/q_{n_k+1})}
{F (p'_{n_k}/q'_{n_k}) - F (p_{n_k}/q_{n_k})}=
\frac{F (p_{n_k}/q_{n_k} \oplus p'_{n_k}/q'_{n_k}) - F (p_{n_k}/q_{n_k})}
{F (p'_{n_k}/q'_{n_k}) - F (p_{n_k}/q_{n_k})}= \\=
\begin{cases}
\frac{c}{\lambda},\quad\text{if } S(p_{n_k}/q_{n_k}) < S(p'_{n_k}/q'_{n_k}). \\
\frac{c}{\lambda^2},\quad\text{if } S(p_{n_k}/q_{n_k}) > S(p'_{n_k}/q'_{n_k}).
\end{cases}
\end{multline}

Then there are two possibilities:

a) $x \in (p_{n_k}/q_{n_k}, (p_{n_k}/q_{n_k} \oplus p'_{n_k}/q'_{n_k})^l )$,

b) $x \in ((p_{n_k}/q_{n_k} \oplus p'_{n_k}/q'_{n_k})^l, p_{n_k}/q_{n_k} \oplus p'_{n_k}/q'_{n_k} )$,

where
$$(p_{n_k}/q_{n_k} \oplus p'_{n_k}/q'_{n_k})^l = p_{n_k}/q_{n_k} \oplus
(p_{n_k}/q_{n_k} \oplus p'_{n_k}/q'_{n_k}),$$
$$(p_{n_k}/q_{n_k} \oplus p'_{n_k}/q'_{n_k})^r = (p_{n_k}/q_{n_k} \oplus
p'_{n_k}/q'_{n_k}) \oplus p'_{n_k}/q'_{n_k}.$$

As the last partial quotient of $p_{n_k}/q_{n_k} \oplus p'_{n_k}/q'_{n_k}$ is equal to $2$,
then by Remark 1, the last partial quotients of $(p_{n_k}/q_{n_k} \oplus p'_{n_k}/q'_{n_k})^r$ and
$(p_{n_k}/q_{n_k} \oplus p'_{n_k}/q'_{n_k})^l$ is not equal to $2$.
As $S(p_{n_k}/q_{n_k}
\oplus p'_{n_k}/q'_{n_k})> S(p_{n_k}/q_{n_k})$,
we deduce from Lemma \ref{ratio} that

\begin{equation}\label{case1a,b}
\frac{F (p'_{n_k+2}/q'_{n_k+2}) - F (p_{n_k+2}/q_{n_k+2})}
{F (p'_{n_k+1}/q'_{n_k+1}) - F (p_{n_k+1}/q_{n_k+1})}=
\begin{cases}
\frac{1}{\lambda},\quad \text{in case a)},\\
\frac{1}{c\lambda}, \quad \text{in case b)}.
\end{cases}
\end{equation}

Consider \textbf{Case 2}. 
Analogously to the case 1 we get
\begin{multline}\label{case_2}
\frac{F (p'_{n_k+1}/q'_{n_k+1}) - F (p_{n_k+1}/q_{n_k+1})}
{F (p'_{n_k}/q'_{n_k}) - F (p_{n_k}/q_{n_k})}=
\frac{F (p'_{n_k}/q'_{n_k}) - F (p_{n_k}/q_{n_k} \oplus p'_{n_k}/q'_{n_k})}
{F (p'_{n_k}/q'_{n_k}) - F (p_{n_k}/q_{n_k})}= \\=
\begin{cases}
\frac{c}{\lambda^2},\quad\text{if } S(p_{n_k}/q_{n_k}) < S(p'_{n_k}/q'_{n_k}). \\
\frac{c}{\lambda},\quad\text{if } S(p_{n_k}/q_{n_k}) > S(p'_{n_k}/q'_{n_k}).
\end{cases}
\end{multline}

We shall consider following 
 subcases:

a) $x \in (p_n/q_n \oplus p'_n/q'_n, (p_n/q_n \oplus p'_n/q'_n)^r )$,

b) $x \in ((p_n/q_n \oplus p'_n/q'_n )^r, p'_n/q'_n )$.

As $S(p_{n_k}/q_{n_k} \oplus p'_{n_k}/q'_{n_k})>S(p'_{n_k}/q'_{n_k})$,
from Lemma \ref{ratio} we see that
\begin{equation}\label{case2a,b}
\frac{F (p'_{n_k+2}/q'_{n_k+2}) - F (p_{n_k+2}/q_{n_k+2})}
{F (p'_{n_k+1}/q'_{n_k+1}) - F (p_{n_k+1}/q_{n_k+1})}=
\begin{cases}
\frac{1}{c\lambda},\quad \text{in case a)}, \\
\frac{1}{\lambda}, \quad \text{in case b)}.
\end{cases}
\end{equation}

As the sequence $\{p_{n_k}/q_{n_k}, p'_{n_k}/q'_{n_k}\}$
is infinite then at least one case
(from the cases 1a),1b),2a),2b)) will occur infinitely often.
So there exists a subsequence $\{p_{n_{k_m}}/q_{n_{k_m}},
p'_{n_{k_m}}/q'_{n_{k_m}}\}$ such that
$$
\frac{F (p'_{n_{k_m}+1}/q'_{n_{k_m}+1}) - F (p_{n_{k_m}+1}/q_{n_{k_m}+1})}
{F (p'_{n_{k_m}}/q'_{n_{k_m}}) - F (p_{n_{k_m}}/q_{n_{k_m}})}= \alpha,
$$
$$
\frac{F (p'_{n_{k_m}+2}/q'_{n_{k_m}+2}) - F (p_{n_{k_m}+2}/q_{n_{k_m}+2})}
{F (p'_{n_{k_m}+1}/q'_{n_{k_m}+1}) - F (p_{n_{k_m}+1}/q_{n_{k_m}+1})}= \beta,
$$
where $\alpha$~--- is one of the numbers $\frac{c}{\lambda}$, $\frac{c}{\lambda^2}$,
and $\beta$~--- is one of the numbers $\frac{1}{\lambda}$, $\frac{1}{c\lambda}$.
Now
$$
G_{n_{k_m}}(x) = \alpha \frac{q_{n_{k_m}+1}q'_{n_{k_m}+1}}{q_{n_{k_m}} q'_{n_{k_m}}},
\quad G_{n_{k_m}+1}(x) = \beta \frac{q_{n_{k_m}+2}q'_{n_{k_m}+2}}{q_{n_{k_m}+1}q'_{n_{k_m}+1}}.
$$
From  (\ref{G_n_k}) we see that
\begin{equation}\label{lim}
\lim_{m \to \infty}
\frac{q_{n_{k_m}} q'_{n_{k_m}}}{q_{n_{k_m}+1}q'_{n_{k_m}+1}}={\alpha},
\quad \lim_{m \to \infty}
\frac{q_{n_{k_m}+1}q'_{n_{k_m}+1}}{q_{n_{k_m}+2}q'_{n_{k_m}+2}}
={\beta}.
\end{equation}
Now we must show that (\ref{lim}) is not possible.
To do this 
%With that end in view 
we distinguish the cases again.

\textbf{1},a) In this case
\begin{gather*}
\left\{p_{{n_k}_m+1}/q_{{n_k}_m+1}, p'_{{n_k}_m+1}/q'_{{n_k}_m+1}\right\} = \left\{
p_{{n_k}_m}/q_{{n_k}_m}, (p_{{n_k}_m}+p'_{{n_k}_m})/(q_{{n_k}_m}+q'_{{n_k}_m})\right\},\\
\left\{p_{{n_k}_m+2}/q_{{n_k}_m+2}, p'_{{n_k}_m+2}/q'_{{n_k}_m+2}\right\} = \left\{
p_{{n_k}_m}/q_{{n_k}_m}, (2p_{{n_k}_m}+p'_{{n_k}_m})/(2q_{{n_k}_m}+q'_{{n_k}_m})\right\}.
\end{gather*}
Now (\ref{lim}) leads to
$$
\lim_{m \to \infty} \frac{q'_{n_{k_m}}}{q_{n_{k_m}} + q'_{n_{k_m}}} = \alpha, \quad
\lim_{m \to \infty} \frac{q_{n_{k_m}} + q'_{n_{k_m}}}{2q_{n_{k_m}} + q'_{n_{k_m}}} = \beta,
$$
where by (\ref{case_1}) and (\ref{case1a,b}) one has $\beta = \frac{1}{\lambda}$,
$\alpha = \frac{c}{\lambda}$ for
$S(p_{{n_k}_m}/q_{{n_k}_m}) < S(p'_{{n_k}_m}/q'_{{n_k}_m})$, and $\alpha = \frac{c}{\lambda^2}$
for $S(p_{{n_k}_m}/q_{{n_k}_m}) > S(p'_{{n_k}_m}/q'_{{n_k}_m})$.
Note that
$$
\frac{2q_{n_{k_m}} + q'_{n_{k_m}}}{q_{n_{k_m}} + q'_{n_{k_m}}} =
2  - \frac{q'_{n_{k_m}}}{q_{n_{k_m}} + q'_{n_{k_m}}}.
$$
So we have
$$
\frac{1}{\beta} =2- \alpha.
$$
For $\beta = \frac{1}{\lambda}$, $\alpha = \frac{c}{\lambda}$ we get
$$
%\lambda = 2 - \frac{c}{\lambda} \Rightarrow \lambda = 2 - \frac{\lambda}
%{1 + \lambda}
%\Rightarrow
\lambda^2 -2 = 0.
$$
For $\beta = \frac{1}{\lambda}$, $\alpha = \frac{c}{\lambda^2}$ we get
$$
%\lambda = 2 - \frac{c}{\lambda^2} \Rightarrow \lambda = 2-\frac{1}{\lambda+1}
%\Rightarrow
\lambda^2 - \lambda - 1 = 0.
$$

In both cases we have a contradiction with the fact that $\lambda$ is a root of equation
(\ref{haract}).

\textbf{1},b) In this case
\begin{equation*}
\left\{p_{{n_k}_m+1}/q_{{n_k}_m+1}, p'_{{n_k}_m+1}/q'_{{n_k}_m+1}\right\} = \left\{
p_{{n_k}_m}/q_{{n_k}_m}, (p_{{n_k}_m}+p'_{{n_k}_m})/(q_{{n_k}_m}+q'_{{n_k}_m})\right\},\\
\end{equation*}
\vspace{-30pt}
\begin{multline*}
\left\{p_{{n_k}_m+2}/q_{{n_k}_m+2}, p'_{{n_k}_m+2}/q'_{{n_k}_m+2}\right\} = \\= \left\{
(2p_{{n_k}_m}+p'_{{n_k}_m})/(2q_{{n_k}_m}+q'_{{n_k}_m}),
(p_{{n_k}_m}+p'_{{n_k}_m})/(q_{{n_k}_m}+q'_{{n_k}_m})\right\},
\end{multline*}
Now (\ref{lim}) leads to
$$
\lim_{m \to \infty} \frac{q'_{n_{k_m}}}{q_{n_{k_m}} + q'_{n_{k_m}}} = \alpha, \quad
\lim_{m \to \infty} \frac{q_{n_{k_m}}}{2q_{n_{k_m}} + q'_{n_{k_m}}} = \beta,
$$
where by (\ref{case_1}) and (\ref{case1a,b}) one  has
$\beta = \frac{1}{c\lambda}$, $\alpha = \frac{c}{\lambda}$ for
$S(p_{{n_k}_m}/q_{{n_k}_m}) < S(p'_{{n_k}_m}/q'_{{n_k}_m})$, and
$\alpha = \frac{c}{\lambda^2}$
for $S(p_{{n_k}_m}/q_{{n_k}_m}) > S(p'_{{n_k}_m}/q'_{{n_k}_m})$.
Note that
$$
\frac{2q_{n_{k_m}} + q'_{n_{k_m}}}{q_{n_{k_m}}} = 1 + \frac{q_{n_{k_m}} + q'_{n_{k_m}}}{q_{n_{k_m}}} = 1 +
\frac{1}{1 - \frac{q'_{n_{k_m}}}{q_{n_{k_m}} + q'_{n_{k_m}}}}.
$$
So we have
$$
\frac{1}{\beta} = 1 + \frac{1}{1-\alpha}
$$
For $\beta = \frac{1}{c\lambda}$, $\alpha = \frac{c}{\lambda}$ we get
$$
%c\lambda = 1 +\frac{1}{1 - \frac{c}{\lambda}} \Rightarrow
%\frac{\lambda}{\lambda - 1} =1+ \frac{1}{1-\frac{\lambda}{\lambda + 1}} \Rightarrow
\lambda^2 - 2 = 0.
$$
For $\beta = \frac{1}{c\lambda}$, $\alpha = \frac{c}{\lambda^2}$
$$
%c\lambda = 1 +\frac{1}{1 - \frac{c}{\lambda^2}} \Rightarrow
%\frac{\lambda}{\lambda - 1} =1+ \frac{1}{1-\frac{1}{\lambda + 1}} \Rightarrow
\lambda^2 - \lambda - 1 = 0.
$$

Again in both cases we have the contradiction with the fact that $\lambda$ is a root of equation
(\ref{haract}).

\textbf{2},a) In this case
$$\left\{p_{{n_k}_m+1}/q_{{n_k}_m+1}, p'_{{n_k}_m+1}/q'_{{n_k}_m+1}\right\} =
\left\{
(p_{{n_k}_m}+p'_{{n_k}_m})/(q_{{n_k}_m}+q'_{{n_k}_m}), p'_{{n_k}_m}/q'_{{n_k}_m}\right\},\\
$$
\vspace{-30pt}
\begin{multline*}
\left\{p_{{n_k}_m+2}/q_{{n_k}_m+2}, p'_{{n_k}_m+2}/q'_{{n_k}_m+2}\right\} =
\\ = \left\{
(p_{{n_k}_m}+p'_{{n_k}_m})/(q_{{n_k}_m}+q'_{{n_k}_m}),
(p_{{n_k}_m}+2p'_{{n_k}_m})/(q_{{n_k}_m}+2q'_{{n_k}_m})\right\}.
\end{multline*}
So (\ref{lim}) leads to
$$
\lim_{m \to \infty} \frac{q_{n_{k_m}}}{q_{n_{k_m}} + q'_{n_{k_m}}} = \alpha, \quad
\lim_{m \to \infty} \frac{q'_{n_{k_m}}}{q_{n_{k_m}} + 2q'_{n_{k_m}}} = \beta,
$$
where by (\ref{case_2}) and (\ref{case2a,b}) one has
$\beta = \frac{1}{c\lambda}$, $\alpha = \frac{c}{\lambda^2}$ for
$S(p_{{n_k}_m}/q_{{n_k}_m}) < S(p'_{{n_k}_m}/q'_{{n_k}_m})$, and
$\alpha = \frac{c}{\lambda}$
for $S(p_{{n_k}_m}/q_{{n_k}_m}) > S(p'_{{n_k}_m}/q'_{{n_k}_m})$.
Note that
$$
\frac{q_{n_{k_m}} + 2q'_{n_{k_m}}}{q'_{n_{k_m}}}  = 1 + \frac{q_{n_{k_m}} + q'_{n_{k_m}}}{q'_{n_{k_m}}} =
1 + \frac{1}{1 - \frac{q_{n_{k_m}}}{q_{n_{k_m}} + q'_{n_{k_m}}}}.
$$
So we have
$$
\frac{1}{\beta} = 1 + \frac{1}{1-\alpha},
$$
and tis case  is reduced to the case 1,b).
%�� ���������������� ����� ����������� ���
%$\beta = \frac{1}{c\lambda}$, $\alpha = \frac{c}{\lambda^2}$ �
%$\beta = \frac{1}{c\lambda}$, $\alpha = \frac{c}{\lambda}$ ��� ���� ��������
%� ������ 1,b).

\textbf{2},b) In this case
\begin{gather*}
\left\{p_{{n_k}_m+1}/q_{{n_k}_m+1}, p'_{{n_k}_m+1}/q'_{{n_k}_m+1}\right\} = \left\{
(p_{{n_k}_m}+p'_{{n_k}_m})/(q_{{n_k}_m}+q'_{{n_k}_m}), p'_{{n_k}_m}/q'_{{n_k}_m}\right\},\\
\left\{p_{{n_k}_m+2}/q_{{n_k}_m+2}, p'_{{n_k}_m+2}/q'_{{n_k}_m+2}\right\} = \left\{
(p_{{n_k}_m}+2p'_{{n_k}_m})/(q_{{n_k}_m}+2q'_{{n_k}_m}),
p'_{{n_k}_m}/q'_{{n_k}_m}\right\}.
\end{gather*}
Now (\ref{lim}) leads to
$$
\lim_{m \to \infty} \frac{q_{n_{k_m}}}{q_{n_{k_m}} + q'_{n_{k_m}}} = \alpha, \quad
\lim_{m \to \infty} \frac{q_{n_{k_m}} + q'_{n_{k_m}}}{q_{n_{k_m}} + 2q'_{n_{k_m}}} = \beta,
$$
where by (\ref{case_2}) and (\ref{case2a,b}) one has
$\beta = \frac{1}{\lambda}$, $\alpha = \frac{c}{\lambda^2}$ for
$S(p_{{n_k}_m}/q_{{n_k}_m}) < S(p'_{{n_k}_m}/q'_{{n_k}_m})$, and
$\alpha = \frac{c}{\lambda}$
for $S(p_{{n_k}_m}/q_{{n_k}_m}) > S(p'_{{n_k}_m}/q'_{{n_k}_m})$.
Note that
$$
\frac{q_{n_{k_m}} + 2q'_{n_{k_m}}}{q_{n_{k_m}} + q'_{n_{k_m}}} =
2 - \frac{q_{n_{k_m}}}{q_{n_{k_m}} + q'_{n_{k_m}}}.
$$
So we have
$$
\frac{1}{\beta} =2- \alpha,
$$
and we have the same situation as in the case 1,a).
%�� ���������������� ����� ����������� ���
%$\beta = \frac{1}{\lambda}$, $\alpha = \frac{c}{\lambda^2}$ �
%$\beta = \frac{1}{\lambda}$, $\alpha = \frac{c}{\lambda}$ ��� ���� ��������
%� ������ 1,a).

Theorem 2 is proved.
\hfill $\blacksquare$

\end{document}